\documentclass[10pt]{article}
\usepackage{amsfonts}
\usepackage{hyperref}
\usepackage{mathrsfs,enumerate,amssymb,amsmath,color,lineno}
\usepackage{indentfirst}
\usepackage{amsmath}
\usepackage{amssymb}
\usepackage[amsmath,thmmarks]{ntheorem}
\usepackage{cite}
\usepackage{graphicx}
\usepackage{tikz}

\newtheorem{definition}{\bf Definition}[section]
\newtheorem{lemma}{\bf Lemma}[section]
\newtheorem{theorem}{\bf Theorem}[section]
\newtheorem{remark}{\bf Remark}[section]
\newtheorem{corollary}{\bf Corollary}[section]
\newtheorem{example}{\bf Example}[section]

\newtheorem{proposition}{\bf Proposition}[section]

\begin{document}
\setcounter{page}{1}

\title{{\textbf{Characterizations of monotone right continuous functions which generate associative functions}}\thanks {Supported by
the National Natural Science Foundation of China (No.12071325)}}
\author{Yun-Mao Zhang\footnote{\emph{E-mail address}: 2115634219@qq.com}, Xue-ping Wang\footnote{Corresponding author. xpwang1@hotmail.com; fax: +86-28-84761502},\\
\emph{School of Mathematical Sciences, Sichuan Normal University,}\\
\emph{Chengdu 610066, Sichuan, People's Republic of China}}

\newcommand{\pp}[2]{\frac{\partial #1}{\partial #2}}
\date{}
\maketitle
\begin{quote}
{\bf Abstract} Associativity of a two-place function $T: [0,1]^2\rightarrow [0,1]$ defined by $T(x,y)=f^{(-1)}(T^*(f(x),f(y)))$ where $T^*:[0,1]^2\rightarrow[0,1]$ is an associative function with neutral element in $[0,1]$, $f: [0,1]\rightarrow [0,1]$ is a monotone right continuous function and $f^{(-1)}:[0,1]\rightarrow[0,1]$ is the pseudo-inverse of $f$ depends only on properties of the range of $f$. The necessary and sufficient conditions for the $T$ to be associative are presented by applying the properties of the monotone right continuous function $f$.

{\textbf{\emph{Keywords}}:} Associative function; Right continuous function; Pseudo-inverse; Semigroup; Triangular norm\\
\end{quote}

\section{Introduction}
Triangular norms (t-norms for short) have been introduced in the framework of statistical  metric spaces for proper generalization of the triangle inequality from metric spaces to statistical metric spaces \cite{CA2006, EP2000, BS1960}. As far as we know, t-norms have played a valuable role in many theoretical and practical applications of mathematics. In the context of fuzzy logic where it is a many-valued propositional logic with a continuum of truth values modelled by unit intervals, t-norms and triangular conorms (t-conorms for short) model the (semantic) interpretations of the logical conjunction and disjunction, respectively. T-norms and t-conorms have also been applied to control, theory of non-additive measures and integrals \cite{EP1995}, measure-free conditioning \cite{UH1999}, decision making \cite{JF1994} and so on.

In the theory of t-norms, the constructing methods of t-norms play an important role. There are many methods in the literature concerning the construction of t-norms, see e.g., \cite{SJ1999, EP2000, AM2004, YO2007, YO2008, PV2005, PV2008, DZ2005}. In particular, based on known t-norms, monotone functions and their pseudo-inverses, Klement, Mesiar and Pap \cite{EP2000} have proved that if a non-decreasing function $f: [0,1]\rightarrow [0,1]$ and a t-norm $T: [0,1]^2\rightarrow [0,1]$ such that
\begin{equation}
\label{eq:1.1}
 T(f(x),f(y))\in \mbox{Ran}(f)\cup [0,f(0^+)]
\end{equation}
for all $x,y \in [0,1)$ and, for all $(x,y)\in [0,1]^2$ with $T(f(x),f(y))\in \mbox{Ran}(f)$,
\begin{equation}
\label{eq:1.2}
 f\circ f^{(-1)}(T(f(x),f(y)))=T(f(x),f(y)).
\end{equation} Then the following function $T_{[f]}: [0,1]^2\rightarrow [0,1]$ given by
 \begin{equation}\label{equ003}
T_{[f]}(x,y)=\left\{
  \begin{array}{ll}
    \min\{x,y\} & \mbox{if }\max\{x,y\}=1, \\
  f^{(-1)}(T(f(x),f(y))) & \mbox{otherwise,}
  \end{array}
\right.
\end{equation}
is a t-norm where $f^{(-1)}$ is the pseudo-inverse of $f$. Furthermore, Klement, Mesiar and Pap \cite{EP2000} used the best-known operations, the usual addition and multiplication of real numbers, to construct t-norms with the help of functions in one variable. They defined an additive generator $f :[0,1]\rightarrow [0,\infty]$ of a t-norm $T$ as a strictly decreasing function which is also right-continuous at $0$ with $f (1) = 0$, such that for all $(x,y)\in [0,1]^2$
\begin{equation}
\label{eq:1.7}
f(x)+f(y)\in \mbox{Ran}(f)\cup [f(0),\infty],
\end{equation}

\begin{equation}
\label{eq:1.8}
 T(x,y)= f^{(-1)}(f(x)+f(y))
\end{equation}
Moreover, they pointed out that the concept of additive generator of t-norms can be further generalized in the sense that a strictly decreasing function $f :[0,1]\rightarrow [0,\infty]$ with $f (1) = 0$ can generate a t-norm $T$ via Eq.\eqref{eq:1.8}, i.e., such that for all $(x,y)\in [0,1]^2$ we have $ T(x,y)= f^{(-1)}(f(x)+f(y))$, without satisfying Eq.\eqref{eq:1.7}. Many results concerning additive generators of associative functions which do not satisfy Eq.\eqref{eq:1.7} can be found \cite{PV1998, PV2005, PV2008, PV2010, PV2013, PV2020, PV2022}. For example, Vicen\'{\i}k \cite{PV2005} showed the necessary and sufficient conditions for associativity of the function $T$ via Eq.\eqref{eq:1.8} in term of properties of the range of an additive generator $f$. Motivated by the work of Vicen\'{\i}k \cite{PV2005} together with Eq.\eqref{equ003}, in this article we consider what is a characterization of all monotone right continuous functions $f: [0,1]\rightarrow [0,1]$ such that the function $T: [0,1]^2\rightarrow [0,1]$ given by \begin{equation}\label{eq5}
T(x,y)= f^{(-1)}(T^*(f(x),f(y)))
\end{equation}
where $f^{(-1)}$ is the pseudo-inverse of $f$ and $T^*: [0,1]^2\rightarrow [0,1]$ is an associative function with neutral element in $[0,1]$ is associative.

The rest of this article is organized as follows. In Section 2, we recall some basic concepts and known results of t-norms and t-conorms, respectively. In Section 3, we give a representation of the range $\mbox{Ran}(f)$ of a non-decreasing right continuous function $f$. In Section 4, we first define an operation $\otimes$ on $\mbox{Ran}(f)$ with $f$ a non-decreasing right continuous function, and then investigate some necessary and sufficient conditions for the operation $\otimes$ being associative. In Section 5, we describe what properties of $\mbox{Ran}(f)$ are equivalent to associativity of the operation $\otimes$. A conclusion is drawn in Section 6.

\section{Preliminaries}

In this section, we recall some known basic concepts and results.
\begin{definition}[\cite{EP2000}]
\emph{A t-norm is a binary operator $T:[0, 1]^2\rightarrow [0, 1]$ such that for all $x, y, z\in[0, 1]$ the following conditions are satisfied:}

$(T1)$  $T(x,y)=T(y,x)$,

$(T2)$  $T(T(x,y),z)=T(x,T(y,z))$,

\emph{$(T3)$  $T(x,y)\leq T(x,z)$ whenever $y\leq z$,}

$(T4)$  $T(x,1)=T(1,x)=x$.
\end{definition}

A binary operator $T:[0, 1]^2\rightarrow [0, 1]$ is called t-subnorm if it satisfies $(T1), (T2), (T3)$, and $T(x,y)\leq \min\{x,y\}$ for all $x,y\in [0, 1]$.

\begin{definition}[\cite{EP2000}]
\emph{A t-conorm is a binary operator $T:[0, 1]^2\rightarrow [0, 1]$ such that for all $x, y, z\in[0, 1]$ the following conditions are satisfied:}

$(S1)$ $S(x,y)=S(y,x)$,

$(S2)$  $S(S(x,y),z)=S(x,S(y,z))$,

\emph{$(S3)$  $S(x,y)\leq S(x,z)$ whenever $y\leq z$,}

$(S4)$ $S(x,0)=S(0,x)=x$.
\end{definition}

\begin{definition}[\cite{EP2000,PV2005}]\label{def2.3}
\emph{Let $p, q, s, t\in [-\infty, \infty]$ with $p<q, s<t$ and $f:[p,q]\rightarrow[s,t]$ be a non-decreasing (resp. non-increasing) function. Then the function $f^{(-1)}:[s,t]\rightarrow[p,q]$ defined by
\begin{equation*}
f^{(-1)}(y)=\sup\{x\in [p,q]\mid f(x)<y\}\,(\mbox{resp. }f^{(-1)}(y)=\sup\{x\in [p,q]\mid f(x)>y\})
\end{equation*}
is called the pseudo-inverse of a non-decreasing (resp. non-increasing) function $f$.}
\end{definition}

Denoted by $A\setminus B=\{x\in A\mid x\notin B\}$ for two sets $A$ and $B$.
\begin{theorem}[\cite{EP2000}]\label{them2.0}
Let $p, q, s, t\in [-\infty, \infty]$ with $p<q, s<t$ and $f:[p,q]\rightarrow[s,t]$ be a non-constant increasing function. Then for each $y\in [s,t]\setminus \emph{Ran}(f)$ we have
\begin{equation*}
\sup\{x\in [p,q]\mid f(x)<y\}=\inf\{x\in [p,q]\mid f(x)>y\},
\end{equation*}
where $\emph{Ran}(f)=\{f(x)\mid x\in[p,q]\}.$
\end{theorem}

\section{The range of a non-decreasing right continuous function}

In this section we give a representation of the range of a non-decreasing right continuous function.

Let $M\subseteq [0,1]$. A point $x\in [0,1]$ is said to be an accumulation point of $M$ from the left if there exists a strictly increasing sequence $\{x_{n}\}_{n\in N}$ of points $x_{n}\in M$ such that $\lim_{n\rightarrow \infty}x_{n}=x$. A point $x\in [0,1]$ is said to be an accumulation point of $M$ from the right if there exists a strictly decreasing sequence $\{x_{n}\}_{n\in N}$ of points $x_{n}\in M$ such that $\lim_{n\rightarrow \infty}x_{n}=x$. A point $x\in [0,1]$ is said to be an accumulation point of $M$ if $x$  is an accumulation point of $M$ from the left or from the right. For a function $f:[0,1]\rightarrow [0,1]$, let $f(a^-)=\lim_ {x \rightarrow a^{- }} f(x)$ for each $a\in (0, 1]$ and $f(a^+)=\lim_ {x \rightarrow a^{+}} f(x)$ for each $a\in [0, 1)$.
Denote
$$\mathcal{A}=\{M \mid \mbox{there is a non-decreasing right continuous function } f:[0,1]\rightarrow[0,1]\mbox{ such that }\mbox{Ran}(f)=M\}.$$
Then we have the following representation of the range of a non-decreasing right continuous function.

\begin{lemma}\label{lem3.1}
Let $M\in \mathcal{A}$ with $M\neq[0, f(1)]$. Then there exist a uniquely determined non-empty countable system $\mathcal{S}=\{[b_k, d_k] \subseteq [0,1]\mid k\in K\}$ of pairwise closed intervals of a positive length which satisfies that for all $[b_k, d_k],[b_l, d_l]\in \mathcal{S}$, $[b_k, d_k]\cap [b_l, d_l]=\emptyset$ or $[b_k, d_k]\cap [b_l, d_l]=\{d_k\}$ when $d_k\leq b_l$, and a uniquely determined non-empty countable set $\mathcal{C}=\{c_k\in [0,1]\mid k\in \overline{K}\}$ such that $[b_k, d_k] \cap \mathcal{C}=\{d_k\}$ or $[b_k, d_k] \cap \mathcal{C}=\{b_k,d_k\}$ for all $k\in K$ and
\begin{equation*}
M= \{c_k\in [0,1]\mid k\in \overline{K}\}\cup \left([0,f(1)]\setminus \left(\bigcup_{k\in K}[b_k, d_k] \right)\right)
\end{equation*}
where $|K|\leq|\overline{K}|$.
\end{lemma}
\begin{proof}We will prove the existence of both $\mathcal{S}$ and $\mathcal{C}$. Since $M\in \mathcal{A}$ with $M\neq[0, f(1)]$, there exists a non-decreasing right continuous function $f:[0,1]\rightarrow[0,1]$
 with $\mbox{Ran}(f)=M$. Take $\mathcal{S}=\{[f(x^-),f(x)]\mid x\in[0,1], f(x^-)<f(x)\}$ and $\mathcal{C}=\{f(x)\mid x\in[0,1], f(x^-)<f(x)\}\cup \{f(x^-)\mid x\in[0,1], f(x^-)<f(x), f(x^-)\in M\}$. It is easy to see that $\mathcal{S}$ is countable. Let $K$ be a countable index set and $|K|=|\mathcal{S}|$. So that we can rewrite $\mathcal{S}$ as $\mathcal{S}=\{[b_k, d_k] \subseteq [0,1]\mid k\in K\}$ in which $b_k=f(k^-)$, $d_k=f(k)$ for each $k\in K$. $\mathcal{C}$ is also countable. Let $\overline{K}$ be a countable index set and $|\overline{K}|=|\mathcal{C}|$. Obviously, $|K|\leq|\overline{K}|$. Rearrange $\mathcal{C}$ from small to large and assign its every element an index $k\in\overline{K}$, for example, the $k$-th element in the rearranged set $\mathcal{C}$ is denoted by $c_k$ with $k\in\overline{K}$. Then $\mathcal{C}=\{c_k\in [0,1]\mid k\in\overline{K}\}$. Notice that for every $k\in K$, there exit $m, n\in\overline{K}$ with $m+1=n$ such that $d_k=f(k)=c_{n}$ and $b_k=f(k^-)=c_{m}$ when $f(k^-)\in M$, and there exits $j\in\overline{K}$ such that $d_k=c_{j}$ when $f(k^-)\notin M$. Obviously, $\mathcal{S}$ and $\mathcal{C}$ have all required properties, respectively.

Now, we will prove the uniqueness of both $\mathcal{S}$ and $\mathcal{C}$. Suppose that a system $\mathcal{S}_1=\{[u_l, v_l] \subseteq [0,1]\mid l\in L\}$ and a set $\mathcal{C}_1=\{c_l\in [0,1]\mid l\in \overline {L}\}$ have all required properties. Next, we will prove that $\mathcal{S}=\mathcal{S}_1$. Fix an arbitrary interval $[b_k, d_k]\in \mathcal{S}$. Choose $t\in [b_k, d_k]$ such that $t\notin M$. Then there exists a $[u_l, v_l]\in\mathcal{S}_1$ such that $t\in [u_l, v_l]$ where $[u_l, v_l] \cap M=\{v_l\}$ or $[u_l, v_l] \cap M=\{u_l,v_l\}$. In what follows, we shall prove that $\mathcal{S}\subseteq \mathcal{S}_1$ by two steps.

Suppose that $u_l<d_k < v_l$. Then it contradicts the facts that $[b_k, d_k] \cap M=\{d_k\}$ or $[b_k, d_k] \cap M=\{b_k,d_k\}$ and $[u_l, v_l] \cap M=\{v_l\}$ or $[u_l, v_l] \cap M=\{u_l,v_l\}$.
If $d_k\leq u_l< v_l$, then it contradicts the facts that $t\in [u_l, v_l]$ and $t\in [b_k, d_k]$ such that $t\notin M$. If $b_k<v_l< d_k $, then it contradicts the facts that $[b_k, d_k] \cap M=\{d_k\}$ or $[b_k, d_k] \cap M=\{b_k,d_k\}$ and $[u_l, v_l] \cap M=\{v_l\}$ or $[u_l, v_l] \cap M=\{u_l,v_l\}$. If $v_l\leq b_k< d_k$, then it contradicts the facts that $t\in [u_l, v_l]$ and $t\in [b_k, d_k]$ such that $t\notin M$. Therefore, $d_k =v_l$.

Suppose that $u_l < b_k < v_l$. Then $b_k=f(x^-)$ for a certain $x\in(0,1]$. If $b_k\in M$, then $b_k\in [u_l, v_l] \cap M$, contrary to $[u_l, v_l] \cap M=\{v_l\}$ or $[u_l, v_l] \cap M=\{u_l,v_l\}$. Thus $b_k\notin M$. Then $b_k$ is an accumulation point of $M$ from the left. This follows that the set $[u_l, v_l] \cap M$ is infinite, contrary to $[u_l, v_l] \cap M=\{v_l\}$ or $[u_l, v_l] \cap M=\{u_l,v_l\}$. If $v_l\leq b_k$, then it contradicts the facts that $t\in [u_l, v_l]$ and $t\in [b_k, d_k]$ such that $t\notin M$. The cases $b_k< u_l<d_k$ and $d_k\leq u_l$ are completely analogous, respectively. Therefore, $b_k=u_l$.

Consequently, $\mathcal{S}\subseteq \mathcal{S}_1$.

Finally, we shall prove that $\mathcal{S}_1\subseteq \mathcal{S}$. Fix an arbitrary interval $[u_l, v_l]\in \mathcal{S}_1$. Choose $t\in [u_l, v_l]$ such that $t\notin M$. Then there exists $k_0\in K$ such that $t\in [b_{k_0}, d_{k_0}]$ with $[b_{k_0}, d_{k_0}]\in \mathcal{S}\subseteq \mathcal{S}_1$. Thus $t\in [u_l, v_l]\cap [b_{k_0}, d_{k_0}]$. Since $[u_k, v_k]\cap [u_l, v_l]=\emptyset$ or $[u_k, v_k]\cap [u_l, v_l]=\{v_k\}$ for all $[u_k, v_k],[u_l, v_l]\in\mathcal{S}_1$, the only possibility is that $[u_l, v_l]=[b_{k_0}, d_{k_0}]\in \mathcal{S}$. Therefore, $\mathcal{S}_1\subseteq \mathcal{S}$.

 Consequently, $\mathcal{S}=\mathcal{S}_1$. This follows that $\mathcal{C}=\mathcal{C}_1$. In particular, both $\mathcal{S}$ and $\mathcal{C}$ are independent of choice of $f$.
\end{proof}

\begin{definition}\label{def3.1}
\emph{Let $M\in \mathcal{A}$. A pair$(\mathcal{S},\mathcal{C})$ is said to be associated with $M\neq [0,f(1)]$ if $\mathcal{S}=\{[b_k, d_k] \subseteq [0,1]\mid k\in K\}$ is a non-empty countable system of pairwise closed intervals of a positive length which satisfies that for all $[b_k, d_k],[b_l, d_l]\in \mathcal{S}$, $[b_k, d_k]\cap [b_l, d_l]=\emptyset$ or $[b_k, d_k]\cap [b_l, d_l]=\{d_k\}$ when $d_k\leq b_l$, and $\mathcal{C}=\{c_k\in [0,1]\mid k\in \overline{K}\}$ is a non-empty countable set such that $[b_k, d_k] \cap \mathcal{C}=\{d_k\}$ or $[b_k, d_k] \cap \mathcal{C}=\{b_k,d_k\}$ for all $k\in K$ and
\begin{equation*}
M= \{c_k\in [0,1]\mid k\in \overline{K}\}\cup \left([0,f(1)]\setminus \left(\bigcup_{k\in K}[b_k, d_k] \right)\right).
\end{equation*}
A pair $(\mathcal{S},\mathcal{C})$ is said to be associated with $M=[0,f(1)]$ if $\mathcal{S}=\{[f(1),f(1)]\}$ and $\mathcal{C}=\{f(1)\}$.}
\end{definition}

We shall briefly write $(\mathcal{S},\mathcal{C})=(\{[b_k, d_k]\mid k\in K\}, \{c_k \mid k\in \overline{K}\})$ instead of $(\mathcal{S},\mathcal{C})=(\{[b_k, d_k] \subseteq [0, 1]\mid k\in K\}, \{c_k\in [0,1] \mid k\in \overline{K}\})$. Denote
$$E=\{c_k\mid \mbox{ there are an } x_0\in[0,1]\mbox{ and }\varepsilon >0\mbox{ such that }f|_{[x_0,x_0+\varepsilon]}=c_k\},$$ and $F=\mathcal{C}\setminus E$, where $f|_{[x_0,x_0+\varepsilon]}$ is a restriction of the function $f$ on $[x_0,x_0+\varepsilon]$.

From Definition \ref{def3.1}, we have the following remark.
\begin{remark} \emph{Let $M\in \mathcal{A}$, $(\mathcal{S},\mathcal{C})=(\{[b_k, d_k]\mid k\in K\}, \{c_k \mid k\in \overline{K}\})$ be associated with $M$ and $f:[0,1]\rightarrow[0,1]$ be a non-decreasing right continuous function with $\mbox{Ran}(f)=M$.
\renewcommand{\labelenumi}{(\roman{enumi})}
\begin{enumerate}
\item For all $a,b\in[0,1]$, if $a<b$ and $b\notin E$, then $[a,b)\cap (M\setminus F)=\emptyset$ if and only if $[a,b]\cap (M\setminus F)=\emptyset$.
\item If $[b_k, d_k]\cap [b_l, d_l]=\{d_k\}$ for some $k,l\in K$, then there exist an $x\in[0,1]$ and $\varepsilon >0$ such that $f\mid_{[x,x+\varepsilon]}=d_k$.
\item If $[b_k, d_k] \cap M=\{b_k,d_k\}$ for a certain $k\in K$, then there exist an $x\in[0,1]$ and $\varepsilon >0$ such that $f\mid_{[x,x+\varepsilon]}=b_k$.
\end{enumerate}}
\end{remark}

\begin{example}\label{exp3.1}\emph{
\renewcommand{\labelenumi}{(\roman{enumi})}
\begin{enumerate}
\item Let the function $g_1:[0,1]\rightarrow [0,1]$ defined by\begin{equation*}
 g_1(x)=\begin{cases}
\frac{1}{2}x & \hbox{if }\ x\in[0,\frac{2}{5}),\\
\frac{2}{5} &  \hbox{if }\ x\in[\frac{2}{5},\frac{3}{5}),\\
x &  \hbox{if }\ x\in[\frac{3}{5},1].\\
\end{cases}
\end{equation*}
Then the pair $(\{[\frac{1}{5},\frac{2}{5}],[\frac{2}{5},\frac{3}{5}]\},\{\frac{2}{5},\frac{3}{5}\})$ is associated with $[0,\frac{1}{5})\cup \{\frac{2}{5}\}\cup [\frac{3}{5},1]\in \mathcal{A}$, $E=\{\frac{2}{5}\}$ and $F=\{\frac{3}{5}\}$.
\item Let the function $g_2:[0,1]\rightarrow [0,1]$ defined by \begin{equation*}
 g_2(x)=\begin{cases}
x & \hbox{if }\ x\in[0,\frac{1}{5}),\\
\frac{1}{5} &  \hbox{if }\ x\in[\frac{1}{5},\frac{2}{5}),\\
\frac{2}{5} &  \hbox{if }\ x\in[\frac{2}{5},\frac{3}{5}),\\
\frac{3}{5} &  \hbox{if }\ x\in[\frac{3}{5},\frac{4}{5}),\\
2(x-\frac{4}{5})+\frac{3}{5} & \hbox{if }\ x\in[\frac{4}{5},1].\\
\end{cases}
\end{equation*}
Then the pair $(\{[\frac{1}{5},\frac{2}{5}],[\frac{2}{5},\frac{3}{5}]\},\{\frac{1}{5},\frac{2}{5},\frac{3}{5}\})$ is associated with $[0,\frac{1}{5}]\cup \{\frac{2}{5}\}\cup [\frac{3}{5},1]\in \mathcal{A}$, $E=\{\frac{1}{5},\frac{2}{5},\frac{3}{5}\}$ and $F=\emptyset$.
\end{enumerate}}
\end{example}

\section{An operation on $\mbox{Ran}(f)$ and its properties}

In this section we first define an operation $\otimes$ on $\mbox{Ran}(f)$ with $f$ a non-decreasing right continuous function, and then investigate some necessary and sufficient conditions for the operation $\otimes$ being associative.

\begin{definition}\label{def3.2}
\emph{Let $M\in \mathcal{A}$. Define a function $G_{M}:[0,1]\rightarrow M $ by
\begin{equation*}
 G_{M}(x)=\begin{cases}
\inf([x,1]\cap M) & \hbox{if }\ x< f(1),\\
f(1) &  \hbox{otherwise}.
\end{cases}
\end{equation*}}
\end{definition}

The next proposition describes the relationships between $M$ and $G_{M}$.

\begin{proposition}\label{prop3.1}
Let $M\in \mathcal{A}$ and $(\mathcal{S},\mathcal{C})=(\{[b_k, d_k]\mid k\in K\}, \{c_k \mid k\in \overline{K}\})$ be associated with $M$. Then for all $x,y\in[0,1]$ and $k\in K$,
\renewcommand{\labelenumi}{(\roman{enumi})}
\begin{enumerate}
\item $G_{M}(x)=x$ if and only if $x\in M$.
\item if $x\notin M$ and $x< f(1)$ then $G_{M}(x)=d_k$ if and only if $x\in [b_k, d_k]$ and $[b_k, d_k] \cap M=\{d_k\}$ or $[b_k, d_k] \cap M=\{b_k,d_k\}$. If $x\notin M$ and $f(1)< x$ then $G_{M}(x)=f(1)$.
\item $G_{M}$ is a non-decreasing function.
\end{enumerate}
\end{proposition}
\begin{proof}
(i) The proof is immediate.

(ii) If $x\notin M$, $x< f(1)$ and $G_{M}(x)=d_k$, then $G_{M}(x)=\inf([x,1]\cap M)$ implies $x\in [b_k, d_k]$ and $d_k\in M$. Therefore, $[b_k, d_k] \cap M=\{d_k\}$ or $[b_k, d_k] \cap M=\{b_k,d_k\}$. The converse implication is obviously available. If $x\notin M$ and $f(1)< x$, then clearly $G_{M}(x)=f(1)$ by Definition \ref{def3.2}.

(iii)
It is immediately by Definition \ref{def3.2}.
\end{proof}

\begin{example}\label{exp3.2}
\emph{In Example \ref{exp3.1},
\renewcommand{\labelenumi}{(\roman{enumi})}
\begin{enumerate}
\item \begin{equation*}
 G_{M}(x)=\begin{cases}
\frac{2}{5} & \hbox{if }\ x\in[\frac{1}{5},\frac{2}{5}),\\
\frac{3}{5} &  \hbox{if }\ x\in(\frac{2}{5},\frac{3}{5}),\\
x &  \mbox{otherwise}.\\
\end{cases}
\end{equation*}
\item \begin{equation*}
 G_{M}(x)=\begin{cases}
\frac{2}{5} & \hbox{if }\ x\in(\frac{1}{5},\frac{2}{5}),\\
\frac{3}{5} &  \hbox{if }\ x\in(\frac{2}{5},\frac{3}{5}),\\
x &  \mbox{otherwise}.\\
\end{cases}
\end{equation*}
\end{enumerate}}
\end{example}

In what follows, we always suppose that $T^*:[0,1]^2\rightarrow [0,1]$ is an associative function with neutral element in $[0,1]$.

\begin{definition}\label{def3.3}
\emph{Let $M\in \mathcal{A}$ and $G_{M}$ be determined by $M$. Define an operation $\otimes:M^2\rightarrow M $ by
\begin{equation*}
x\otimes y=G_{M}(T^*(x,y)).
\end{equation*}}
\end{definition}
\begin{example}\label{exp3.3}
\emph{In Example \ref{exp3.2},
\renewcommand{\labelenumi}{(\roman{enumi})}
\begin{enumerate}
\item \begin{equation*}
 x\otimes y=\begin{cases}
\frac{2}{5} & \hbox{if }\ T^*(x,y)\in[\frac{1}{5},\frac{2}{5}),\\
\frac{3}{5} &  \hbox{if }\ T^*(x,y)\in(\frac{2}{5},\frac{3}{5}),\\
T^*(x,y) &  \hbox{otherwise}.\\
\end{cases}
\end{equation*}
\item \begin{equation*}
 x\otimes y=\begin{cases}
\frac{2}{5} & \hbox{if }\ T^*(x,y)\in(\frac{1}{5},\frac{2}{5}),\\
\frac{3}{5} &  \hbox{if }\ T^*(x,y)\in(\frac{2}{5},\frac{3}{5}),\\
T^*(x,y) &  \hbox{otherwise}.\\
\end{cases}
\end{equation*}
\end{enumerate}}
\end{example}

\begin{proposition}\label{prop3.2}
Let $M\in \mathcal{A}$ and $(\mathcal{S},\mathcal{C})=(\{[b_k, d_k]\mid k\in K\}, \{c_k \mid k\in \overline{K}\})$ be associated with $M$. Then for all $x,y\in M$ and $k\in K$,
\renewcommand{\labelenumi}{(\roman{enumi})}
\begin{enumerate}
\item $x\otimes y=T^*(x,y)$ if and only if $T^*(x,y)\in M$.
\item if $T^*(x,y)\notin M$ and $T^*(x,y)< f(1)$, then $x\otimes y=d_k$ if and only if $T^*(x,y)\in [b_k, d_k]$ and $[b_k, d_k] \cap M=\{d_k\}$ or $[b_k, d_k] \cap M=\{b_k,d_k\}$. If $T^*(x,y)\notin M$ and $f(1)< T^*(x,y)$ then $x\otimes y=f(1)$.
\end{enumerate}
\end{proposition}
\begin{proof}
It is an immediate consequence of Propositions \ref{prop3.1} (i) and (ii) and Definition \ref{def3.3}.
\end{proof}

Let $[a,b]\subseteq [-\infty,\infty]$ with $a\leq b$. Then by convention, $\sup \emptyset=a$ and $\inf \emptyset=b$.
\begin{lemma}\label{lem3.2}
Let $M\in \mathcal{A}$, $(\mathcal{S},\mathcal{C})=(\{[b_k, d_k]\mid k\in K\}, \{c_k \mid k\in \overline{K}\})$ be associated with $M$ and $f:[0,1]\rightarrow [0,1] $ be a non-decreasing right continuous function with $\emph{Ran}(f)=M$. Then $G_{M}(x)=f(f^{(-1)}(x))$ for all $x\in[0,1]$.
\end{lemma}
\begin{proof}
If $x\in \mbox{Ran}(f)$, then $f(f^{(-1)}(x))=x$ since $f$ is right continuous.

If $x\notin \mbox{Ran}(f)$ and $x< f(1)$, then taking $m=\inf\{t\in \mbox{Ran}(f)\mid x<t\}$. From Theorem \ref{them2.0} and $f$ is right continuous, we have
\begin{eqnarray*}
f(f^{(-1)}(x))&=&f(\sup\{y\in [0,1]\mid f(y)<x\})\\
&=&f(\inf\{y\in [0,1]\mid x< f(y)\})\\
&=&\inf \{f(y)\in [0,1]\mid x< f(y)\}\\
&=&m
\end{eqnarray*}
and there exists a $k\in K$ such that $x\in [b_k, d_k]$ with $d_k=\inf\{t\in \mbox{Ran}(f)\mid x<t\}$. From Proposition \ref{prop3.1}, we get $G_{M}(x)=m$. If $f(1)<x$, then $G_{M}(x)=f(1)$ and $f(f^{(-1)}(x))=f(\sup\{y\in [0,1]\mid f(y)<x\})=f(1)$.
Thus, $G_{M}(x)=f(f^{(-1)}(x))$ for all $x\in[0,1]$.
\end{proof}

Note that it is clear that $G_{M}(x)=f(f^{(-1)}(x))$ for all $x\in[0,1]$ when the function $f$ is a constant function. Therefore, we next  always suppose that the function $f$ is a non-constant non-decreasing right continuous function.
\begin{remark}\label{rem3.2}\emph{In Lemma \ref{lem3.2}, the condition that the function $f:[0,1]\rightarrow [0,1] $ is right continuous can not be deleted when $f$ is not strictly monotone.
For instance, let the function $f:[0,1]\rightarrow [0,1]$ be defined by \begin{equation*}
 f(x)=\begin{cases}
\frac{1}{2}x & \hbox{if }\ x\in[0,\frac{2}{5}],\\
\frac{2}{5} &  \hbox{if }\ x\in(\frac{2}{5},\frac{3}{5}],\\
x &  \hbox{if }\ x\in(\frac{3}{5},1].\\
\end{cases}
\end{equation*}
It is clear that $f$ is not right continuous and $G_{M}(\frac{2}{5})=\frac{2}{5}>\frac{1}{5}=f(f^{(-1)}(\frac{2}{5}))$, i.e., $G_{M}\neq f(f^{(-1)})$.}
\end{remark}

Denote
$$H=\{\min \{x\in[0,1]\mid f(x)=y\}\mid y\in E\}, J=\{x\in [0,1]\mid f(x)\in M\setminus E \}\mbox{ and }
D=H\cup J.$$

\begin{definition}\label{def3.4} \emph{Let $f:[0,1]\rightarrow [0,1]$ be a non-decreasing right continuous function.
Define a function $f^*:D\rightarrow [0,1]$ by
\begin{equation*}
f^*(x)=f(x)
\end{equation*}
for all $x\in D$, a function $d^*:[0,1]\rightarrow D$ by
\begin{equation*}
d^*(x)=f^{(-1)}(x)
\end{equation*}
 for all $x\in [0,1]$ and a two-place function $T_0:D^2\rightarrow D$ by
 \begin{equation*}
T_0(x,y)=d^*(T^*(f^*(x),f^*(y)))
\end{equation*}}
\end{definition}
for all $x,y\in D$, respectively.

\begin{proposition}\label{prop3.3}
Let $f:[0,1]\rightarrow [0,1]$ be a non-decreasing right continuous function. Then $f^*$ is a strictly increasing function, $d^*(f^*(x))=x$ for all $x\in D$ and $f^*(d^*(y))=f(f^{(-1)}(y))$ for all $y\in [0,1]$.
\end{proposition}
\begin{proof}
It follows from Definition \ref{def3.4} that $f^*$ is strictly increasing. Now, let $x\in D$. Then from Definition \ref{def2.3}, $d^*(f^*(x))=\sup\{t\in D\mid f(t)<f^*(x)\}=x$. Because of $d^*(y)=f^{(-1)}(y)$ for all $y\in [0,1]$, $f^*(d^*(y))=f(f^{(-1)}(y))$ for all $y\in [0,1]$.
\end{proof}

\begin{proposition}\label{prop3.4}
Let $f:[0,1]\rightarrow [0,1]$ be a non-decreasing right continuous function. Then
 \begin{equation*}
x\otimes y=f^*(T_0(d^*(x),d^*(y)))
\end{equation*}
for all $x,y\in M$, and
 \begin{equation*}
T_0(x,y)=d^*(f^*(x)\otimes f^*(y))
\end{equation*}
for all $x,y\in D$. Moreover, $T_0$ is associative if and only if $\otimes$ is associative.
\end{proposition}
\begin{proof} For each $x,y\in M$, from Definitions \ref{def3.3} and \ref{def3.4}, Lemma \ref{lem3.2} and Proposition \ref{prop3.3}, we have
\begin{eqnarray*}
x\otimes y&=&G_{M}(T^*(x,y))\\
&=&f(f^{(-1)}(T^*(x,y)))\\
&=&f^*(d^*(T^*(x,y)))\\
&=&f^*(d^*(T^*(f^*(d^*(x)),f^*(d^*(y)))))\\
&=&f^*(T_0(d^*(x),d^*(y))).
\end{eqnarray*}
Because of $x,y\in M$, there exist $u,v\in D$ such that $f^*(u)=x$, $f^*(v)=y$. Thus, from Proposition \ref{prop3.3}, we have
\begin{equation*}
f^*(u)\otimes f^*(v)=x\otimes y=f^*(T_0(d^*(x),d^*(y)))=f^*(T_0(u,v)).
\end{equation*}
This follows that $d^*(f^*(u)\otimes f^*(v))=d^*(f^*(T_0(u,v))=T_0(u,v)$.

It is obviously that $T_0$ is associative if and only if $\otimes$ is associative.
\end{proof}

\begin{proposition}\label{prop3.5}Let $f:[0,1]\rightarrow [0,1]$ be a non-decreasing right continuous function and $T:[0,1]^2\rightarrow [0,1]$ be a function defined by Eq.(\ref{eq5}).
Then $T_0$ is associative if and only if $T$ is associative.
\end{proposition}
\begin{proof}
Suppose that $T_0$ is associative. Let us prove that $T(T(x,y),z)=T(x,T(y,z))$ for all $x,y,z\in[0,1]$. Note that from Definition \ref{def3.4} and Proposition \ref{prop3.3}, if $x,y\in D$ then \begin{eqnarray*}
T(x,y)&=&f^{(-1)}(T^*(f(x),f(y)))\\
&=&d^*(T^*(f^*(x),f^*(y)))\\
&=&T_0(x,y).
\end{eqnarray*}
If $x\notin D$ and $y\in D$ then there exists a $u\in D$ such that $f^*(u)=f(x)$ and
\begin{eqnarray*}
T(x,y)&=&f^{(-1)}(T^*(f(x),f(y)))\\
&=&f^{(-1)}(T^*(f^*(u),f(y)))\\
&=&d^*(T^*(f^*(u),f^*(y)))\\
&=&T_0(u,y).
\end{eqnarray*}
If $x\in D$ and $y\notin D$ then there exists a $v\in D$ such that $f^*(v)=f(y)$ and
\begin{eqnarray*}
T(x,y)&=&f^{(-1)}(T^*(f(x),f(y)))\\
&=&f^{(-1)}(T^*(f^*(x),f(v)))\\
&=&d^*(T^*(f^*(x),f^*(v)))\\
&=&T_0(x,v).
\end{eqnarray*}
If $x\notin D$ and $y\notin D$ then there exist two elements $u,v\in D$ such that $f^*(u)=f(x)$ and $f^*(v)=f(y)$, respectively, and
\begin{eqnarray*}
T(x,y)&=&f^{(-1)}(T^*(f(x),f(y)))\\
&=&f^{(-1)}(T^*(f^*(u),f^*(v)))\\
&=&d^*(T^*(f^*(u),f^*(v)))\\
&=&T_0(u,v).
\end{eqnarray*}

In the following, for completing the proof of associativity, we distinguish four cases.
\begin{itemize}
  \item Let $x,y,z\in [0,1]$ and $x,y,z\in D$. Then
\begin{eqnarray*}
T(T(x,y),z)&=&T(T_0(x,y),z)\\
&=&T_0(T_0(x,y),z)\\
&=&T_0(x,T_0(y,z))\\
&=&T(x,T_0(y,z))\\
&=&T(x,T(y,z)).
\end{eqnarray*}
\item Let $x,y,z\in [0,1]$ and two coordinates of the triple $(x,y,z)$ be exactly contained in $D$. There are three subcases. If $x\notin D$, $y,z\in D$, then there exists $u\in D$ such that $f^*(u)=f(x)$ and
\begin{eqnarray*}
T(T(x,y),z)&=&T(T_0(u,y),z)\\
&=&T_0(T_0(u,y),z)\\
&=&T_0(u,T_0(y,z))\\
&=&T_0(u,T(y,z))\\
&=&T(x,T(y,z)).
\end{eqnarray*}
The other two subcases $y\notin D$, $x,z\in D$ and $z\notin D$, $x,y\in D$ are completely analogous.
 \item Let $x,y,z\in [0,1]$ and exact one element of $\{x,y,z\}$ be contained in $D$. If $x,y\notin D$, $z\in D$, then there exist two elements $u,v\in D$ such that $f^*(u)=f(x)$ and $f^*(v)=f(y)$, respectively, and
\begin{eqnarray*}
T(T(x,y),z)&=&T(T_0(u,v),z)\\
&=&T_0(T_0(u,v),z)\\
&=&T_0(u,T_0(v,z))\\
&=&T_0(u,T(y,z))\\
&=&T(x,T(y,z)).
\end{eqnarray*}
 The other two subcases $y,z\notin D$, $x\in D$ and $x,z\notin D$, $y\in D$ are completely analogous.
 \item Let $x,y,z\in [0,1]$ and $x,y,z\notin D$. Then there exist three elements $u,v,s\in D$ such that $f^*(u)=f(x)$, $f^*(v)=f(y)$ and $f^*(s)=f(z)$, respectively, and \begin{eqnarray*}
T(T(x,y),z)&=&T(T_0(u,v),z)\\
&=&T_0(T_0(u,v),s)\\
&=&T_0(u,T_0(v,s))\\
&=&T_0(u,T(y,z))\\
&=&T(x,T(y,z)).
\end{eqnarray*}
\end{itemize}

Conversely, it is obviously true from Definition \ref{def3.4}.
\end{proof}

The following corollary is a simple consequence of Propositions \ref{prop3.4} and \ref{prop3.5}.
\begin{corollary}\label{coro3.1}
Let $f:[0,1]\rightarrow [0,1]$ be a non-decreasing right continuous function and $T:[0,1]^2\rightarrow [0,1]$ be a function defined by Eq.(\ref{eq5}).
Then $T$ is associative if and only if $\otimes$ is associative.
\end{corollary}

\section{Associativity of the operation $\otimes$}

In this section we show a necessary and sufficient condition for associativity of the operation $\otimes$ which answers what properties of $M$ are equivalent to associativity of the operation $\otimes$.

Let $M\subseteq [0,1]$. Define $O(M) = \bigcup_{x,y\in M}[\min\{x,y\}, \max\{x,y\})$ when $M\neq \emptyset$ (where $[x, x)=\emptyset$), and $O(M)=\emptyset$ when $M=\emptyset$. It is clear that $O(M)=\emptyset$ if and only if $M=\emptyset$ or $M$ contains only one element.
Let $\emptyset\neq A,B\subseteq [0,1]$ and $c\in[0,1]$. Denote $T^*(A,c)=\{T^*(x,c)\mid x\in A\}$, $T^*(c,A)=\{T^*(c,x)\mid x\in A\}$ and $T^*(A,B)=\{T^*(x,y)\mid x\in A, y\in B\}$.

Note that if $f(1)=1$, then $T^*(x,y)\leq 1=f(1)$ for all $x,y\in[0,1]$. If $f(1)<1$, then there may exist $x,y\in[0,1]$ such that $T^*(x,y)>f(1)$. Therefore, in what follows, we continue our discussion by distinguishing two cases $f(1)=1$ and $f(1)<1$.

We first consider the case $f(1)=1$. The following definition is needed.
\begin{definition}\label{def3.5}
\emph{Let $M\in \mathcal{A}$ with $f(1)=1$ and $(\mathcal{S},\mathcal{C})=(\{[b_k, d_k]\mid k\in K\}, \{c_k \mid k\in \overline{K}\})$ be associated with $M$. For all $y\in M$, and $k,l\in K$, define $$M_{k}^{y}=\{x\in M \mid \mbox{ either }T^*(x,y)\in (b_k, d_k)\mbox{ when } b_k\in M\mbox{ or }T^*(x,y)\in [b_k, d_k)\mbox{ when } b_k\notin M \},$$
$$M_{y}^{k}=\{x\in M \mid \mbox{ either }T^*(y,x)\in (b_k, d_k)\mbox{ when } b_k\in M\mbox{ or }T^*(y,x)\in [b_k, d_k)\mbox{ when } b_k\notin M \},$$
  $$M^{y}=\{x\in M \mid T^*(x,y)\in M\},$$ $$M_{y}=\{x\in M \mid T^*(y,x)\in M\},$$ $$I_{k}^{y}=\{d_k\}\cup T^*(M_{k}^{y},y),$$ $$I_{y}^{k}=\{d_k\}\cup T^*(y,M_{y}^{k})\mbox{ and}$$
\begin{equation*}
J_{k,l}^{y}=\begin{cases}
O(T^*(M_{k}^{y},d_l)\cup T^*(d_k,M_{y}^{l})), & \hbox{if }\ M_{k}^{y}\neq\emptyset,  M_{y}^{l}\neq\emptyset,\\
\emptyset, & \hbox{otherwise.}
\end{cases}
\end{equation*}
Put $$\mathfrak{T_{1}}(M)=\bigcup_{y\in M}\bigcup_{k\in K}\bigcup_{t\in M_{y}}O(T^*(I_{k}^{y},t)),$$ $$\mathfrak{T_{2}}(M)=\bigcup_{y\in M}\bigcup_{k\in K}\bigcup_{t\in M^{y}}O(T^*(t,I_{y}^{k})),$$ $$\mathfrak{T_{3}}(M)=\bigcup_{y\in M}\bigcup_{k,l\in K}J_{k,l}^{y}\mbox{ and}$$ $$\mathfrak{T}(M)=\mathfrak{T_{1}}(M)\cup\mathfrak{T_{2}}(M)\cup\mathfrak{T_{3}}(M).$$}
\end{definition}

In the following, we further suppose that $T^*:[0,1]^2\rightarrow [0,1]$ is monotone. Then we have the following two lemmas.
\begin{lemma}\label{lem3.3}
Let $M\in \mathcal{A}$ with $f(1)=1$ and $(\mathcal{S},\mathcal{C})=(\{[b_k, d_k]\mid k\in K\}, \{c_k \mid k\in \overline{K}\})$ be associated with $M$. Let $M_1, M_2\subseteq [0,1]$ be two non-empty sets and $c\in[0,1]$.
Then\\
(1) $O(T^*(M_1\cup M_2,c)\cap M\neq\emptyset$ if and only if there exist $x\in M_1$ and $y\in M_2$ such that $$[\min\{T^*(x,c),T^*(y,c)\},\max\{T^*(x,c),T^*(y,c)\})\cap M \neq\emptyset.$$
(2) $O(T^*(c,M_1\cup M_2)\cap M\neq\emptyset$ if and only if there exist $x\in M_1$ and $y\in M_2$ such that $$[\min\{T^*(c,x),T^*(c,y)\},\max\{T^*(c,x),T^*(c,y)\})\cap M\neq\emptyset.$$
\end{lemma}
\begin{proof}
We just prove (1), the proof of statement (2) being analogous. Suppose that $T^*(O(M_1\cup M_2),c)\cap M\neq\emptyset$. Note that if $T^*$ is non-decreasing, then there exist two elements $s,t\in M_1\cup M_2$ with $s<t$ such that $[T^*(s,c),T^*(t,c))\cap M\neq\emptyset$. If $T^*$ is non-increasing, then there exist two elements $s,t\in M_1\cup M_2$ with $s<t$ such that $[T^*(t,c),T^*(s,c))\cap M\neq\emptyset$.

The following proof is split into three cases.
\renewcommand{\labelenumi}{(\roman{enumi})}
\begin{enumerate}
\item The case when one of $\{s,t\}$ is contained in $M_1$ and the other is contained in $M_2$ is obvious.

\item The case that $s,t\in M_1$ with $s<t$. Choose $w\in M_2$. There are three subcases when $T^*$ is non-decreasing as below.
\renewcommand{\labelenumi}{(\roman{enumi})}
\begin{enumerate}
\item If $w\leq s$ then put $x=t$ and $y =w$.
\item If $w\geq t$ then put $x=s$ and $y =w$.
\item If $s<w<t$, then $[T^*(s,c),T^*(w,c))\cap M\neq\emptyset$ or $[T^*(w,c),T^*(t,c))\cap M\neq\emptyset$. If $[T^*(s,c),T^*(w,c))$ $\cap M\neq\emptyset$ then put $x=s$ and $y =w$. If $[T^*(w,c),T^*(t,c))\cap M\neq\emptyset$ then put $x=t$ and $y =w$.
\end{enumerate}
In any subcase when $T^*$ is non-decreasing, one can check that $$[\min\{T^*(x,c),T^*(y,c)\},\max\{T^*(x,c),T^*(y,c)\})\cap M \neq\emptyset.$$ There are three subcases when $T^*$ is non-increasing as below.
\begin{enumerate}
\item If $w\leq s$ then put $x=t$ and $y =w$.
\item If $w\geq t$ then put $x=s$ and $y =w$.
\item If $s<w<t$, then $[T^*(t,c),T^*(w,c))\cap M\neq\emptyset$ or  $[T^*(w,c),T^*(s,c))\cap M\neq\emptyset$. If $[T^*(t,c),T^*(w,c))$ $\cap M\neq\emptyset$ then put $x=t$ and $y =w$. If $[T^*(w,c),T^*(s,c))\cap M\neq\emptyset$ then put $x=s$ and $y =w$.
\end{enumerate}
In any subcase when $T^*$ is non-increasing, one can check that $$[\min\{T^*(x,c),T^*(y,c)\},\max\{T^*(x,c),T^*(y,c)\})\cap M \neq\emptyset.$$
\item The case $s,t\in M_2$ is completely analogous to (ii).
\end{enumerate}

The converse implication is obvious.
\end{proof}

\begin{lemma}\label{lem3.5}
Let $M\in \mathcal{A}$ with $f(1)=1$ and $(\mathcal{S},\mathcal{C})=(\{[b_k, d_k]\mid k\in K\}, \{c_k \mid k\in \overline{K}\})$ be associated with $M$. Then for any $x,y\in[0,1]$, $G_M(x)\neq G_M(y)$ if and only if $[\min\{x,y\},\max\{x,y\})\cap M\neq\emptyset$.
\end{lemma}
\begin{proof}
If $G_M(x)\neq G_M(y)$, then from Proposition \ref{prop3.1}, $G_M(\min\{x,y\})=s<t=G_M(\max\{x,y\})$ with $s,t\in M$. Obviously $[s,t)\cap M\neq\emptyset$. Choose $w\in [s,t)\cap M$. From Proposition \ref{prop3.1} and $w=G_M(w)$, we have $G_M(\min\{x,y\})\leq G_M(w)<G_M(\max\{x,y\})$, thus, $\min\{x,y\}\leq w<\max\{x,y\}$. Hence $[\min\{x,y\},\max\{x,y\})\cap M\neq\emptyset$.

Conversely, supposing $[\min\{x,y\},\max\{x,y\})\cap M\neq\emptyset$, then $x\neq y$, say $x<y$. We need to consider the following four situations. If $x\in M, y\in M$, then from Proposition \ref{prop3.1}, $G_M(x)<G_M(y)$. If $x\in M, y\notin M$, then there exists a $k\in K$ such that $y\in[b_k,d_k]$ and $G_M(y)=d_k>y$. Thus from Proposition \ref{prop3.1}, $G_M(x)<y<d_k=G_M(y)$. If $x\notin M, y\in M$, then there exists a $k\in K$ such that $x\in[b_k,d_k]$. Thus from Proposition \ref{prop3.1}, $G_M(x)=d_k\leq y=G_M(y)$. We further claim that $G_M(x)=d_k< y=G_M(y)$. Otherwise, $d_k=y$ implies $[x,y)\cap M=\emptyset$, a contradiction. If $x\notin M, y\notin M$, then there exist two elements $k,l \in K$ with $k<l$ such that $x\in[b_k,d_k]$ and $y\in[b_l,d_l]$. Therefore, $G_M(x)=d_k<d_l=G_M(y)$ by Proposition \ref{prop3.1}. In summary, we always have $G_M(x)< G_M(y)$.
\end{proof}

The following theorem characterizes what properties of $M$ with $f(1)=1$ are equivalent to associativity of the operation $\otimes$.
\begin{theorem}\label{them3.1}
Let $M\in \mathcal{A}$ with $f(1)=1$ and $(\mathcal{S},\mathcal{C})=(\{[b_k, d_k]\mid k\in K\}, \{c_k \mid k\in \overline{K}\})$ be associated with $M$. Then $(M,\otimes)$ is a semigroup if and only if $\mathfrak{T}(M)\cap M=\emptyset$.
\end{theorem}
\begin{proof} Let $(\mathcal{S},\mathcal{C})=(\{[b_k, d_k]\mid k\in K\}, \{c_k \mid k\in \overline{K}\})$ be associated with $M\in \mathcal{A}$. In order to complete the proof, it is enough to prove that the operation $\otimes$ on $M$ is not associative if and only if $\mathfrak{T}(M)\cap M\neq\emptyset$.

 Suppose that the operation $\otimes$ is not associative, i.e., there exist three elements $x,y,z\in M$ such that $(x\otimes y)\otimes z\neq x\otimes(y\otimes z)$. Then we claim that $T^*(x,y)\notin M$ or $T^*(y,z)\notin M$. Otherwise, from Definition \ref{def3.3}, $T^*(x,y)\in M$ and $T^*(y,z)\in M$ would imply $(x\otimes y)\otimes z=G_M(T^*(T^*(x,y),z))= G_M(T^*(x,T^*(y,z)))= x\otimes(y\otimes z)$, a contradiction. The following proof is split into three cases.

(i) Let $T^*(x,y)\notin M$ and $T^*(y,z)\in M$. Then $y\otimes z=T^*(y,z)$ and there exists a $k\in K$ such that $T^*(x,y)\in (b_k,d_k)$ or $T^*(x,y)\in [b_k,d_k)$ with $b_k\notin M$. Thus $x\otimes y=d_k$ and from Lemma \ref{lem3.1}, $M\cap [b_k,d_k]=\{d_k\}$ or $M\cap [b_k,d_k]=\{b_k,d_k\}$. It follows from Definition \ref{def3.3} that $G_M(T^*(d_k,z))=(x\otimes y)\otimes z\neq x\otimes(y\otimes z)=G_M(T^*(x,T^*(y,z)))$. On the other hand, by associativity of $T^*$, we have $G_M(T^*(x,T^*(y,z)))=G_M(T^*(T^*(x,y),z))$. Thus $G_M(T^*(d_k,z))\neq G_M(T^*(T^*(x,y),z))$. Therefore, by Lemma \ref{lem3.5}, $[\min\{T^*(d_k,z),T^*(T^*(x,y),z)\},\max\{T^*(d_k,z),T^*(T^*(x,y),z)\})\cap M\neq\emptyset $. Obviously, $\bigcup_{t\in M_{y}}O(T^*(I_{k}^{y},t))=\bigcup_{t\in M_{y}}\bigcup_{m,n\in T^*(\{d_k\}\cup T^*(M_{k}^{y},y),t)}[\min\{m,n\}, \max\{m,n\})$, and $x\in M_{k}^{y}$. So that $[\min\{T^*(d_k,z),T^*(T^*(x,y),z)\},\max\{T^*(d_k,z),T^*(T^*(x,y),z)\})\subseteq \bigcup_{t\in M_{y}}O(T^*(I_{k}^{y},t))$, which implies $\bigcup_{t\in M_{y}}O(T^*(I_{k}^{y},t))\cap M \neq\emptyset$.

(ii) Let $T^*(x,y)\in M$ and $T^*(y,z)\notin M$. In completely analogous to (i), $\bigcup_{t\in M_{y}}O(T^*(t,I_{y}^{k}))\cap M \neq\emptyset$.
(iii) Let $T^*(x,y)\notin M$ and $T^*(y,z)\notin M$. Then from $T^*(x,y)\notin M$, there exists a $k\in K$ such that $T^*(x,y)\in (b_k,d_k)$ or $T^*(x,y)\in [b_k,d_k)$ with $b_k\notin M$. Thus $x\otimes y=d_k$ and from Lemma \ref{lem3.1}, $M\cap [b_k,d_k]=\{d_k\}$ or $M\cap [b_k,d_k]=\{b_k,d_k\}$. From $T^*(y,z)\notin M$, there exists an $l\in K$ such that $T^*(y,z)\in (b_l,d_l)$ or $T^*(y,z)\in [b_l,d_l)$ with $b_l\notin M$. Hence $y\otimes z=d_l$ and from Lemma \ref{lem3.1}, $M\cap [b_l,d_l]=\{d_l\}$ or $M\cap [b_l,d_l]=\{b_l,d_l\}$. Therefore, $G_M(T^*(d_k,z))=(x\otimes y)\otimes z\neq x\otimes(y\otimes z)=G_M(T^*(x,d_l))$. Applying Lemma \ref{lem3.5}, $[\min\{T^*(d_k,z),T^*(x,d_l)\}, \max\{T^*(d_k,z),T^*(x,d_l)\})\cap M\neq \emptyset$. Obviously, $x\in M_{k}^{y}$ and $T^*(x,d_l)\in T^*(M_{k}^{y},d_l)$. Similarly, $z\in M_{y}^{l}$ and $T^*(d_k,z)\in T^*(d_k,M_{y}^{l})$. Therefore, $[\min\{T^*(d_k,z),T^*(x,d_l)\},\max\{T^*(d_k,z),T^*(x,d_l)\})\subseteq J_{k,l}^{y}$. This follows $J_{k,l}^{y}\cap M\neq \emptyset$.

From (i), (ii) and (iii), we deduce that $\mathfrak{T}(M)\cap M\neq\emptyset$.

Conversely, suppose $\mathfrak{T}(M)\cap M\neq\emptyset$. Then there exist a $y\in M$ and two elements $k,l\in K$ such that $\bigcup_{t\in M_{y}}O(T^*(I_{k}^{y},t))\cap M \neq\emptyset$ or $\bigcup_{t\in M_{y}}O(T^*(t,I_{y}^{k}))\cap M \neq\emptyset$ or $J_{k,l}^{y}\cap M\neq \emptyset$. We distinguish three cases as follows.

(i) If $\bigcup_{t\in M_{y}}O(T^*(I_{k}^{y},t))\cap M \neq\emptyset$, then there exists a $z\in M_{y}$ such that $O(T^*(I_{k}^{y},z))\cap M \neq\emptyset$. Applying Lemma \ref{lem3.3}, there exist $u\in\{d_k\}$, $v\in T^*(M_{k}^{y},y)$ such that
$$[\min\{T^*(u,z),T^*(v,z)\},\max\{T^*(u,z),T^*(v,z)\})\cap M\neq\emptyset.$$
Because of $v\in T^*(M_{k}^{y},y)$, there exists an $x\in M_{k}^{y}$ such that $T^*(x,y)=v$. Therefore, there exist $u\in\{d_k\}$ and $x\in M_{k}^{y}$ such that
$$[\min\{T^*(d_k,z),T^*(T^*(x,y),z)\},\max\{T^*(d_k,z),T^*(T^*(x,y),z)\})\cap M\neq\emptyset.$$
Consequently, from Lemma \ref{lem3.5} we have $G_M(T^*(d_k,z))\neq G_M(T^*(T^*(x,y),z))$.

On the other hand, from $x\in M_{k}^{y}$ we have $T^*(x,y)\in (b_k,d_k)$ or $T^*(x,y)\in [b_k,d_k)$ with $b_k\notin M$. Thus $x\otimes y=d_k$ and from Lemma \ref{lem3.1}, $M\cap [b_k,d_k]=\{d_k\}$ or $M\cap [b_k,d_k]=\{b_k,d_k\}$. From $z\in M_{y}$, we have $T^*(y,z)\in M$. This follows $y\otimes z=T^*(y,z)$.

Therefore, $(x\otimes y)\otimes z=G_M(T^*(d_k,z))\neq G_M(T^*(T^*(x,y),z))=G_M(T^*(x,T^*(y,z)))= x\otimes(y\otimes z)$.

(ii) If $\bigcup_{t\in M_{y}}O(T^*(t,I_{y}^{k}))\cap M \neq\emptyset$, then in complete analogy to (i), $(x\otimes y)\otimes z\neq x\otimes(y\otimes z)$.

(iii) If $J_{k,l}^{y}\cap M\neq \emptyset$, then $J_{k,l}^{y}\neq \emptyset$. By the definition of $J_{k,l}^{y}$, $T^*(O(T^*(M_{k}^{y},d_l)\cup T^*(d_k,M_{y}^{l})),e)\cap M\neq \emptyset$ where $e$ is a neutral element of $T^*$, $T^*(M_{k}^{y},d_l)\neq\emptyset$ and $T^*(d_k,M_{y}^{l})\neq\emptyset$. Applying Lemma \ref{lem3.3}, there exist $u\in T^*(d_k,M_{y}^{l})$, $v\in T^*(M_{k}^{y},d_l)$ such that $$[\min\{T^*(u,e),T^*(v,e)\},\max\{T^*(u,e),T^*(v,e)\})\cap M\neq\emptyset.$$
Because $u\in T^*(d_k,M_{y}^{l})$ and $v\in T^*(M_{k}^{y},d_l)$, there exist an $x\in M_{k}^{y}$ and a $z\in M_{y}^{l}$ such that $u=T^*(d_k,z)$, $v=T^*(x,d_l)$. Therefore, there exist an $x\in M_{k}^{y}$ and a $z\in M_{y}^{l}$ such that
$$[\min\{T^*(d_k,z),T^*(x,d_l)\},\max\{T^*(d_k,z),T^*(x,d_l)\})\cap M\neq\emptyset$$
since $e$ is a neutral element of $T^*$. Further, by Lemma \ref{lem3.5} we have $G_M(T^*(d_k,z))\neq G_M(T^*(x,d_l))$.

On the other hand, from $x\in M_{k}^{y}$ we have $T^*(x,y)\in (b_k,d_k)$ or $T^*(x,y)\in [b_k,d_k)$ with $b_k\notin M$. Thus $x\otimes y=d_k$ and from Lemma \ref{lem3.1}, $M\cap [b_k,d_k]=\{d_k\}$ or $M\cap [b_k,d_k]=\{b_k,d_k\}$.
From $z\in M_{y}^{l}$ we have $T^*(y,z)\in (b_l,d_l)$ or $T^*(y,z)\in [b_l,d_l)$ with $b_l\notin M$. Thus $y\otimes z=d_l$ and from Lemma \ref{lem3.1}, $M\cap [b_l,d_l]=\{d_l\}$ or $M\cap [b_l,d_l]=\{b_l,d_l\}$.

Therefore, $(x\otimes y)\otimes z=G_M(T^*(d_k,z))\neq G_M(T^*(x,d_l))= x\otimes(y\otimes z)$.
\end{proof}

\begin{corollary}\label{coro3.2}
Let $M\in \mathcal{A}$ with $f(1)=1$, $(\mathcal{S},\mathcal{C})=(\{[b_k, d_k]\mid k\in K\}, \{c_k \mid k\in \overline{K}\})$ be associated with $M$ and $T^*$ be commutative. Then $(M,\otimes)$ is a semigroup if and only if $(\mathfrak{T_{1}}(M)\cup\mathfrak{T_{3}}(M)) \cap M=\emptyset$.
\end{corollary}
\begin{proof}It is easy to see that $M_{k}^{y}=M_{y}^{k}$, $M_{y}=M^{y}$, $I_{k}^{y}=I_{y}^{k}$ and $\mathfrak{T_{1}}(M)=\mathfrak{T_{2}}(M)$ if $T^*$ is commutative. Therefore, from Theorem \ref{them3.1}, $(M,\otimes)$ is a semigroup if and only if $(\mathfrak{T_{1}}(M)\cup\mathfrak{T_{3}}(M)) \cap M=\emptyset$.
\end{proof}

Further, from Corollary \ref{coro3.1} and Theorem \ref{them3.1} we have the following corollary.

\begin{corollary}\label{coro3.3}
Let $f:[0,1]\rightarrow [0,1]$ be a non-decreasing right continuous function with $f(1)=1$ and $T:[0,1]^2\rightarrow [0,1]$ be a function defined by Eq.(\ref{eq5}). Then the function $T$ is associative if and only if $\mathfrak{T}(M)\cap M=\emptyset$.
\end{corollary}

Next, we consider the case $f(1)<1$. We need the following definition.

\begin{definition}\label{def3.6}
\emph{Let $M\in \mathcal{A}$ with $f(1)<1$ and $(\mathcal{S},\mathcal{C})=(\{[b_k, d_k]\mid k\in K\}, \{c_k \mid k\in \overline {K}\})$ be associated with $M$. For all $y\in M$, and $k,l\in K$, set $K'=K\cup\{\tau\}$ where $\tau\notin K$. Define $$\mathfrak{M}_{k}^{y}=\{x\in M \mid \mbox{ either }T^*(x,y)\in (b_k, d_k)\mbox{ when } b_k\in M\mbox{ or }T^*(x,y)\in [b_k, d_k)\mbox{ when } b_k\notin M\},$$ $$\mathfrak{M}_{\tau}^{y}=\{x\in M \mid f(1)\leq T^*(x,y)\},$$
$$\mathfrak{M}_{y}^{k}=\{x\in M \mid \mbox{ either }T^*(y,x)\in (b_k, d_k)\mbox{ when } b_k\in M\mbox{ or }T^*(y,x)\in [b_k, d_k)\mbox{ when } b_k\notin M\},$$ $$\mathfrak{M}_{y}^{\tau}=\{x\in M \mid f(1)\leq T^*(y,x)\},$$
  $$\mathfrak{M}^{y}=\{x\in M \mid T^*(x,y)\in M\setminus \{f(1)\}\},$$ $$\mathfrak{M}_{y}=\{x\in M \mid T^*(y,x)\in M\setminus \{f(1)\}\},$$
  $$\mathfrak{I}_{k}^{y}=\{d_k\}\cup T^*(\mathfrak{M}_{k}^{y},y),$$ $$\mathfrak{I}_{\tau}^{y}=\{f(1)\}\cup T^*(\mathfrak{M}_{\tau}^{y},y),$$
  $$\mathfrak{I}_{y}^{k}=\{d_k\}\cup T^*(y,\mathfrak{M}_{y}^{k}),$$ $$\mathfrak{I}_{y}^{\tau}=\{f(1)\}\cup T^*(y,\mathfrak{M}_{y}^{\tau}),$$  \begin{equation*}
\mathfrak{J}_{k,l}^{y}=\begin{cases}
O(T^*(\mathfrak{M}_{k}^{y},d_l)\cup T^*(d_k,\mathfrak{M}_{y}^{l})), & \hbox{if }\ \mathfrak{M}_{k}^{y}\neq\emptyset,  \mathfrak{M}_{y}^{l}\neq\emptyset,\\
\emptyset, & \hbox{otherwise.}
\end{cases}
\end{equation*}
\begin{equation*}
\mathfrak{J}_{\tau,l}^{y}=\begin{cases}
O(T^*(\mathfrak{M}_{\tau}^{y},d_l)\cup T^*(f(1),\mathfrak{M}_{y}^{l})), & \hbox{if }\ \mathfrak{M}_{\tau}^{y}\neq\emptyset,  \mathfrak{M}_{y}^{l}\neq\emptyset,\\
\emptyset, & \hbox{otherwise.}
\end{cases}
\end{equation*}
\begin{equation*}
\mathfrak{J}_{k,\tau}^{y}=\begin{cases}
O(T^*(\mathfrak{M}_{k}^{y},f(1))\cup T^*(d_k,\mathfrak{M}_{y}^{\tau})), & \hbox{if }\ \mathfrak{M}_{k}^{y}\neq\emptyset,  \mathfrak{M}_{y}^{\tau}\neq\emptyset,\\
\emptyset, & \hbox{otherwise.}
\end{cases}
\end{equation*}
\begin{equation*}
\mathfrak{J}_{\tau,\tau}^{y}=\begin{cases}
O(T^*(\mathfrak{M}_{\tau}^{y},f(1))\cup T^*(f(1),\mathfrak{M}_{y}^{\tau})), & \hbox{if }\ \mathfrak{M}_{\tau}^{y}\neq\emptyset,  \mathfrak{M}_{y}^{\tau}\neq\emptyset,\\
\emptyset, & \hbox{otherwise.}
\end{cases}
\end{equation*} $$I_{1}(M)=\bigcup_{y\in M}\bigcup_{k\in K'}\bigcup_{t\in \mathfrak{M}_{y}}O(T^*(\mathfrak{I}_{k}^{y},t)),$$ $$I_{2}(M)=\bigcup_{y\in M}\bigcup_{k\in K'}\bigcup_{t\in \mathfrak{M}^{y}}O(T^*(t,\mathfrak{I}_{y}^{k})),$$ $$I_{3}(M)=\bigcup_{y\in M}\bigcup_{k,l\in K'}\mathfrak{J}_{k,l}^{y}\mbox{ and}$$ $$I(M)=I_{1}(M)\cup I_{2}(M)\cup I_{3}(M).$$}
\end{definition}

In complete analogy to Lemma \ref{lem3.3}, we can prove the following lemma.
\begin{lemma}\label{lem3.6}
Let $M\in \mathcal{A}$ with $f(1)<1$ and $(\mathcal{S},\mathcal{C})=(\{[b_k, d_k]\mid k\in K\}, \{c_k \mid k\in \overline{K}\})$ be associated with $M$. Let $M_1, M_2\subseteq [0,1]$ be two non-empty sets and $c\in[0,1]$.
Then\\
(1) $O(T^*(M_1\cup M_2,c)\cap (M\setminus \{f(1)\})\neq\emptyset$ if and only if there exist $x\in M_1$ and $y\in M_2$ such that $$[\min\{T^*(x,c),T^*(y,c)\},\max\{T^*(x,c),T^*(y,c)\})\cap (M\setminus \{f(1)\}) \neq\emptyset.$$
(2) $O(T^*(c,M_1\cup M_2))\cap (M\setminus \{f(1)\})\neq\emptyset$ if and only if there exist $x\in M_1$ and $y\in M_2$ such that $$[\min\{T^*(c,x),T^*(c,y)\},\max\{T^*(c,x),T^*(c,y)\})\cap (M\setminus \{f(1)\})\neq\emptyset.$$
\end{lemma}
\begin{lemma}\label{lem3.7}
Let $M\in \mathcal{A}$ with $f(1)<1$ and $(\mathcal{S},\mathcal{C})=(\{[b_k, d_k]\mid k\in K\}, \{c_k \mid k\in \overline{K}\})$ be associated with $M$. Then for any $x,y\in[0,1]$, $G_M(x)\neq G_M(y)$ if and only if $[\min\{x,y\},\max\{x,y\})\cap (M\setminus \{f(1)\})\neq\emptyset$.
\end{lemma}
\begin{proof}
If $G_M(x)\neq G_M(y)$, then from Proposition \ref{prop3.1}, $G_M(\min\{x,y\})=s<t=G_M(\max\{x,y\})$ with $s,t\in (M\setminus \{f(1)\})$. Obviously $[s,t)\cap (M\setminus \{f(1)\})\neq\emptyset$. Choose $w\in [s,t)\cap (M\setminus \{f(1)\})$. From Proposition \ref{prop3.1} and $w=G_M(w)$, we have $G_M(\min\{x,y\})\leq G_M(w)<G_M(\max\{x,y\})$, thus, $\min\{x,y\}\leq w<\max\{x,y\}$. Hence $[\min\{x,y\},\max\{x,y\})\cap (M\setminus \{f(1)\})\neq\emptyset$.

Conversely, supposing $[\min\{x,y\},\max\{x,y\})\cap (M\setminus \{f(1)\})\neq\emptyset$, then $x\neq y$, say $x<y$. We need to consider the following four situations. In the case $x\in M, y\in M$, from Proposition \ref{prop3.1}, $G_M(x)<G_M(y)$. In the case $x\in M, y\notin M$, if $x<f(1)$, $y<f(1)$, then there exists a $k\in K$ such that $y\in[b_k,d_k]$ and $G_M(y)=d_k>y$. Thus from Proposition \ref{prop3.1}, $G_M(x)<y<d_k=G_M(y)$. If $x<f(1)$, $y>f(1)$, then $G_M(x)=x<f(1)=G_M(y)$. If $x=f(1)$, $y>f(1)$, then $[x,y)\cap (M\setminus \{f(1)\})=\emptyset$, a contradiction. In the case $x\notin M, y\in M$, there exists a $k\in K$ such that $x\in[b_k,d_k]$. Thus from Proposition \ref{prop3.1}, $G_M(x)=d_k\leq y=G_M(y)$. We further claim that $G_M(x)=d_k< y=G_M(y)$. Otherwise, $d_k=y$ implies $[x,y)\cap (M\setminus \{f(1)\})=\emptyset$, a contradiction. In the case $x\notin M, y\notin M$, if $x<f(1)$, $y<f(1)$, then there exist two elements $k,l \in K$ with $k<l$ such that $x\in[b_k,d_k]$ and $y\in[b_l,d_l]$. Therefore, $G_M(x)=d_k<d_l=G_M(y)$ by Proposition \ref{prop3.1}. If $x<f(1)$, $y>f(1)$, then we claim that there exists a $z\in M\setminus \{f(1)\}$ such that $z=\min\{t\in (M\setminus \{f(1)\})| t>x\}$. Otherwise, $[x,y)\cap (M\setminus \{f(1)\})=\emptyset$, a contradiction. Therefore, $G_M(x)=z<f(1)=G_M(y)$. In summary, we always have $G_M(x)< G_M(y)$.
\end{proof}

The following theorem characterizes what properties of $M$ with $f(1)<1$ are equivalent to associativity of the operation $\otimes$.
\begin{theorem}\label{them3.2}
Let $M\in \mathcal{A}$ with $f(1)<1$ and $(\mathcal{S},\mathcal{C})=(\{[b_k, d_k]\mid k\in K\}, \{c_k \mid k\in \overline{K}\})$ be associated with $M$. Then $(M,\otimes)$ is a semigroup if and only if $I(M)\cap (M\setminus \{f(1)\})=\emptyset$.
\end{theorem}
\begin{proof}This is shown in complete analogy to the proof of Theorem \ref{them3.1}.
\end{proof}

From Corollary \ref{coro3.1} and Theorem \ref{them3.2} we have the following corollary.
\begin{corollary}\label{coro3.4}
Let $f:[0,1]\rightarrow [0,1]$ be a non-decreasing right continuous function with $f(1)<1$ and $T:[0,1]^2\rightarrow [0,1]$ be a function defined by Eq.(\ref{eq5}). Then the function $T$ is associative if and only if $I(M)\cap (M\setminus \{f(1)\})=\emptyset$.
\end{corollary}

Note that, when $T^*$ is a t-norm, we always have $T^*(f(x),f(y))\leq \min \{f(x),f(y)\}\leq f(1)$, and in this case, Theorem \ref{them3.1} is suitable for all non-decreasing right continuous functions $f$. Therefore, we have the following corollary.
\begin{corollary}\label{coro3.5}
Let $f:[0,1]\rightarrow [0,1]$ be a non-decreasing right continuous function and $T:[0,1]^2\rightarrow [0,1]$ be a function defined by Eq.(\ref{eq5}). Then if $T^*$ is a t-norm, then $T$ is a t-subnorm if and only if $\mathfrak{T}(M)\cap M=\emptyset$.
\end{corollary}

Further, from Corollary \ref{coro3.5} we can get the following corollary.
\begin{corollary}\label{coro3.6}
Let $M\in \mathcal{A}$, $(\mathcal{S},\mathcal{C})=(\{[b_k, d_k]\mid k\in K\}, \{c_k \mid k\in \overline{K}\})$ be associated with $M$ and $T^*$ be t-norm. If $(M,\otimes)$ is a semigroup and $T^*(M,M)\subseteq M\cup [0,f(0)]$, then $\mathfrak{T}(M)\subseteq [0,f(0)]$ where $f:[0,1]\rightarrow [0,1]$ is a non-decreasing right continuous function.
\end{corollary}

In particular, we have the following remark.
\begin{remark}\label{rem3.3} \emph{Let $f:[0,1]\rightarrow [0,1]$ be a non-decreasing right continuous function and $T:[0,1]^2\rightarrow [0,1]$ be a function defined by Eq.(\ref{eq5}).}
\begin{enumerate}\renewcommand{\labelenumi}{(\roman{enumi})}
\item \emph{ If $T^*$ is a t-norm, then $T(x,1)=f^{(-1)}(T^*(f(x),f(1)))\leq f^{(-1)}(f(x))\leq x$ for all $x\in [0,1]$. So, $1$ is not necessary a neutral element of $T$. Therefore, if $\mathfrak{T}(M)\cap M=\emptyset$, then $T$ is not necessary a t-norm. However, if $f$ is further a strictly increasing function with $f(1)=1$, then $T$ is a t-norm. Another way is to slightly modify the function $T$ as for all $x,y\in [0,1]$, \begin{equation}\label{eq:1.3}
T(x,y)=\left\{
  \begin{array}{ll}
    \min\{x,y\} & \hbox{if }\max\{x,y\}=1, \\
  f^{(-1)}(T^*(f(x),f(y))) & \hbox{otherwise.}
  \end{array}
\right.
\end{equation} Then one can check that $T$ is a t-norm.}

\item \emph{If $T^*$ is a t-conorm, then $T(x,0)=f^{(-1)}(T^*(f(x),f(0)))\geq f^{(-1)}(f(x))$ and $f^{(-1)}(f(x))\leq x$ for all $x\in [0,1]$. So, $0$ is not necessary a neutral element of $T$. Therefore, if $f(1)=1$ and $\mathfrak{T}(M)\cap M=\emptyset$, or $f(1)<1$ and $I(M)\cap (M\setminus \{f(1)\})=\emptyset$, then $T$ is not necessary a t-conorm. However, if $f$ is further a strictly increasing function with $f(0)=0$, then $T$ is a t-conorm.}
\end{enumerate}
\end{remark}

\begin{example}\label{exap5.1} \emph{Let $M\in \mathcal{A}$ and $(\mathcal{S},\mathcal{C})=(\{[b_k, d_k]\mid k\in K\}, \{c_k \mid k\in \overline{K}\})$ be associated with $M$.}
\renewcommand{\labelenumi}{(\roman{enumi})}
\begin{enumerate}
\item \emph{Let $T^*$ be a commutative, non-decreasing and associative binary function with neutral element $e\in[0,1]$, and the function $f:[0,1]\rightarrow [0,1]$ be defined by \begin{equation*}
 f(x)=\begin{cases}
0 & \hbox{if }\ x\in[0,a),\\
1 &  \hbox{otherwise}.
\end{cases}
\end{equation*}
Then $M=\{0,1\}$, $\mathfrak{T_{1}}(M)=\mathfrak{T_{2}}(M)=\mathfrak{T_{3}}(M)=\emptyset$. So, $\mathfrak{T}(M) \cap M=\emptyset$ and $(M,\otimes)$ is a semigroup by Theorem \ref{them3.1}. Moreover, by Eq.~\eqref{eq:1.3} $T=T_D$ (a drastic t-norm) when $a=1$.}

\item \emph{Let $T^*$ be a t-norm, fix elements $a,c\in(0,1)$ and the right continuous function $f:[0,1]\rightarrow [0,1]$ be defined by \begin{equation*}
 f(x)=\begin{cases}
\frac{c}{a}x & \hbox{if }\ x<a,\\
1 &  \hbox{otherwise}.
\end{cases}
\end{equation*}
Then $M=[0,c)\cup \{1\}$, $\mathfrak{T_{1}}(M)=\mathfrak{T_{2}}(M)=\mathfrak{T_{3}}(M)=\emptyset$. So, $\mathfrak{T}(M) \cap M=\emptyset$ and $(M,\otimes)$ is a semigroup by Theorem \ref{them3.1}. Moreover, from Remark \ref{rem3.3} (i),
\begin{equation*}
T(x,y)=\left\{
  \begin{array}{ll}
    \frac{a}{c}\cdot T^*(\frac{c}{a}x,\frac{c}{a}y) & \hbox{if }(x,y)\in[0,a)^{2}, \\
    a & \hbox{if }(x,y)\in[a,1)^{2}, \\
    \min\{x,y\} & \hbox{otherwise.}
  \end{array}
\right.
\end{equation*}}
\item \emph{Let $T^*(x,y)=xy$ for all $x,y\in[0,1]$, the right continuous function $f:[0,1]\rightarrow [0,1]$ be strictly increasing and $\mbox{Ran}(f)=M$. If $M=\{0\}\cup(\bigcup_{n\in N}[e^{1-2n},e^{2-2n}))\cup \{1\}$, then $\mathfrak{T_{1}}(M)=\mathfrak{T_{2}}(M)=\mathfrak{T_{3}}(M)=\bigcup_{n\in N}[e^{-2n},e^{1-2n})$ and $\mathfrak{T}(M)\cap M=\emptyset$. Therefore, $(M,\otimes)$ is a semigroup by Theorem \ref{them3.1}.}
\item \emph{Let $T^*(x,y)=x+y-xy$ for all $x,y\in[0,1]$, and the function $f:[0,1]\rightarrow [0,1]$ be defined by \begin{equation*}
 f(x)=\begin{cases}
0 & \hbox{if }\ x\in[0,1),\\
\frac{1}{2} &  \hbox{otherwise}.
\end{cases}
\end{equation*}
Then $M=\{0,\frac{1}{2}\}$. Obvious, $T^*(\frac{1}{2},\frac{1}{2})=\frac{3}{4}>\frac{1}{2}=f(1)$, $I_1(M)=I_2(M)=\emptyset$ and $I_3(M)=[\frac{1}{2},\frac{3}{4})$. So, $I(M) \cap (M\setminus\{f(1)\})=\emptyset$ and $(M,\otimes)$ is a semigroup by Theorem \ref{them3.2}. In fact, from Eq.(\ref{eq5}),
\begin{equation*}
T(x,y)=\left\{
  \begin{array}{ll}
    0 & \hbox{if }(x,y)\in[0,1)^{2}, \\
    1 & \hbox{otherwise.}
  \end{array}
\right.
\end{equation*}
It is easy to check that $T$ is an associative function.}
\end{enumerate}
\end{example}
\section{Conclusions}
In this article, we proved that a function $T:[0,1]^2\rightarrow[0,1]$ defined by Eq.(\ref{eq5}) is associative if and only if $\mathfrak{T}(M)\cap M=\emptyset$ (resp. $I(M)\cap (M\setminus \{f(1)\})=\emptyset$) where $M=\mbox{Ran}(f)$ with $f$ a non-decreasing right continuous function and $f(1)=1$ (resp. $f(1)<1$). It is worth pointing out that one can verify that all results obtained are hold when $f$ is a non-increasing right continuous function. Therefore, we answered what are characterizations of all monotone functions $f: [0,1]\rightarrow [0,1]$ such that the function $T: [0,1]^2\rightarrow [0,1]$ given by Eq.(\ref{eq5}) is associative when $f$ is a monotone right continuous function. One may be interesting whether our results are true or not when $f$ is just a monotone function. Up to now, this problem is still open even if we have Remark \ref{rem3.2} in hand. Another interesting question is the relationships of our results to the work of Vicen\'{\i}k \cite{PV2005}. It is evident that our results will be false when we substitute $T^*$ by the usual addition operation ``+" since the usual addition operation ``+" is not a function from $[0,1]^2$ to $[0,1]$. The biggest difference in essence is that our results are true when $f$ is a monotone right continuous function, and the work of Vicen\'{\i}k \cite{PV2005} hold under the condition that $f$ is a strictly monotone function while a monotone right continuous function needs not to be a strictly monotone one as are shown in Examples \ref{exap5.1} (i), (ii) and (iv).
%\section*{Acknowledgments}
%The authors are grateful to the anonymous referees for their valuable comments and suggestions,
%which help the authors to improve the earlier version.

%%%%%%%%%%%%%%%%%%%%%%%%%%%%%%%%%%%%%%%%%%%%%%%%%%%%%%%%%%%%%%%%%%%%%%%%%%%%%%%%%%%%%%%
\section*{Declaration of competing interest}
The authors declare that they have no known competing financial interests or personal relationships that could have appeared to influence the work reported in this paper.
%%%%%%%%%%%%%%%%%%%%%%%%%%%%%%%%%%%%%%%%%%%%%%%%%%%%%%%%%%%%%%%%%%%%%%%%%%%%%%%%%%%%

\end{document}